# NUMERICAL SCHEMES FOR NONLINEAR PREDICTOR FEEDBACK

**Iasson Karafyllis[*] and Miroslav Krstic[**]**

[*]Dept. of Environmental Eng., Technical University of Crete, 73100, Chania, Greece, email: ikarafyl@enveng.tuc.gr

[**]Dept. of Mechanical and Aerospace Eng., University of California, San Diego, La Jolla, CA 92093-0411, U.S.A., email: krstic@ucsd.edu

**Abstract**
Implementation is a common problem with feedback laws with distributed delays. This paper focuses on a specific aspect of the implementation problem for predictor-based feedback laws: the problem of the approximation of the predictor mapping. It is shown that the numerical approximation of the predictor mapping by means of a numerical scheme in conjunction with a hybrid feedback law that uses sampled measurements, can be used for the global stabilization of all forward complete nonlinear systems that are globally asymptotically stabilizable and locally exponentially stabilizable in the delay-free case. Special results are provided for the linear time invariant case. Explicit formulae are provided for the estimation of the parameters of the resulting hybrid control scheme.

## 1. Introduction

Feedback laws with distributed delays arise when predictor-based methodologies are applied to systems with input or measurement delays. The literature for predictor-based feedback laws is reviewed in [12] and the recent works [8,9,11,12,13,14] have extended the predictor-based methodologies to nonlinear systems and time-varying systems. The reader can consult [21,29] for detailed reviews of linear control systems with distributed delays.

   A common problem with feedback laws with distributed delays is implementation. The numerical approximation of distributed delays by discrete delays can lead to instability: this problem was first presented in [26] in the context of predictor-based feedback. The problem was considered as an important open problem in [23] and there are many works which are devoted to the solution of the implementation problem: see [15,16,17,18,19,20,22,24,27,28,29]. Most of the results are focused on the linear case.

   This paper focuses on a different aspect of the implementation problem for predictor-based feedback laws: the problem of the approximation of the predictor mapping. This problem is very common for nonlinear systems: for nonlinear systems it is very rare that the solution map is known analytically. The recent work [5] was devoted to the approximation of the predictor mapping with a successive approximation approach: the method is suitable for globally Lipschitz systems. The problem of approximation of the predictor mapping is a very important aspect of the implementation problem for predictor-based feedback laws but it is different from the usual problem of the approximation of distributed delays by discrete delays. The latter problem will not be studied in the present work.

The idea of the numerical approximation of the predictor mapping by means of a numerical scheme for solving ordinary differential equations arises naturally as a possible method for solving the problem of approximation of the predictor mapping. However, certain obstructions exist, which are not encountered in standard numerical analysis results. The first obstruction is the existence of inputs: control theory tackles systems with inputs (control systems), whereas standard results in numerical analysis are dealing with dynamical systems (systems without inputs). Exception is the work [2] (see also references therein). A second problem is the scarcity of explicit formulae for the approximation error (which coincides with the so-called global discretization error in numerical analysis): in most cases the estimates of the approximation error are qualitative (see [3,4]).

In this work, we show that the numerical approximation of the predictor mapping by means of a numerical scheme for solving ordinary differential equations is indeed one methodology that can be used with success for systems which are not globally Lipschitz. More specifically, we focus on the explicit Euler scheme. We also study the sampling problem, i.e., the problem where the measurement is not available on-line but it is only available at discrete time instants (the sampling times). The problem is solved by means of a hybrid feedback law and the main result is given next.

**Theorem 1.1:** *Consider the delay-free system:*

$$\dot{x}(t) = f(x(t), u(t))$$
$$x(t) \in \Re^n, u(t) \in \Re^m \qquad (1.1)$$

*where* $f : \Re^n \times \Re^m \to \Re^n$ *is a continuously differentiable mapping with* $f(0,0) = 0$. *Assume that:*

**(A1)** *System (1.1) is forward complete.*

**(A2)** *There exists a continuously differentiable function* $k \in C^1(\Re^n; \Re^m)$ *with* $k(0) = 0$ *such that* $0 \in \Re^n$ *is a Globally Asymptotically Stable and Locally Exponentially Stable equilibrium point of the closed-loop system (1.1) with* $u(t) = k(x(t))$.

*Then for every* $\tau > 0$, $r > 0$ *there exists a locally bounded mapping* $N : \Re^n \times \Re_+ \to \{1, 2, 3, ...\}$, *a constant* $\omega > 0$ *and a locally Lipschitz, non-decreasing function* $C : \Re_+ \to \Re_+$ *with* $C(0) = 0$, *such that for every partition* $\{T_i\}_{i=0}^{\infty}$ *of* $\Re_+$ *with* $\sup_{i \geq 0}(T_{i+1} - T_i) \leq r$, *for every* $x_0 \in \Re^n$ *and* $u_0 \in L^{\infty}([-\tau, 0); U)$, *the solution* $(x(t), u(t)) \in \Re^n \times \Re^m$ *of the closed-loop system*

$$\dot{x}(t) = f(x(t), u(t - \tau))$$
$$x(t) \in \Re^n, u(t) \in \Re^m \qquad (1.2)$$

*with*

$$\dot{z}(t) = f(z(t), k(z(t))) \,, \quad z(t) \in \Re^n$$
$$u(t) = k(z(t))$$
, for $t \in [T_i, T_{i+1})$ (1.3)

and

$$z(T_i) = z_N \tag{1.4}$$

where $N := N\left(x(T_i), \sup_{T_i - \tau \leq s < T_i} |u(s)|\right)$, $h = \frac{\tau}{N}$ and

$$z_{i+1} = z_i + \int_{ih}^{(i+1)h} f(z_i, u(T_i - \tau + s)) ds \,, \text{ for } i = 0, ..., N-1 \text{ and } z_0 = x(T_i) \tag{1.5}$$

and initial condition $x(0) = x_0$ and $u(s) = u_0(s)$ for $s \in [-\tau, 0)$ satisfies the following inequality for all $t \geq 0$:

$$|x(t)| + \sup_{t - \tau \leq s < t} |u(s)| \leq \exp(-\omega t) C\left(|x_0| + \sup_{-\tau \leq s < 0} |u_0(s)|\right) \tag{1.6}$$

The notions of Global Asymptotic Stability and Local Exponential stability employed in the statement of Theorem 1.1 are the standard notions used in the literature (see [10]). The notion of forward completeness for (1.1) is the standard notion that guarantees existence of the solution of (1.1) for all times, all initial conditions and all possible inputs (see [1]). Notice that (1.5) is the application of the explicit Euler numerical scheme to the control system (1.2) with step size $h = \frac{\tau}{N}$. Since the number of the grid points $N := N\left(x(T_i), \sup_{T_i - \tau \leq s < T_i} |u(s)|\right)$ is a function of the state and the input, it is clear that different time steps are used each time that a new measurement arrives. Theorem 1.1 is proved by means of a combined Lyapunov and small-gain methodology and its proof is constructive. In Section 3 the control practitioner will find explicit formulae for the computation of $N := N\left(x(T_i), \sup_{T_i - \tau \leq s < T_i} |u(s)|\right)$, which require the knowledge of an appropriate Lyapunov function for the closed-loop system (1.1) with $u(t) = k(x(t))$. Notice that the fact that the function $C: \Re_+ \to \Re_+$ involved in (1.6) is locally Lipschitz with $C(0) = 0$ guarantees the analogue of local exponential stability for complicated systems such as the closed-loop system (1.2) with (1.3), (1.4) and (1.5) (systems with delays and hybrid features), since the estimate

$$|x(t)| + \sup_{t - \tau \leq s < t} |u(s)| \leq \Omega \exp(-\omega t)\left(|x_0| + \sup_{-\tau \leq s < 0} |u_0(s)|\right), \text{ for all } t \geq 0$$

holds for certain appropriate constant $\Omega > 0$ and for initial conditions with $|x_0| + \sup_{-\tau \leq s < 0} |u_0(s)|$ sufficiently small. Therefore, both global asymptotic stability and local exponential stability are preserved, despite the delay, the sampled measurements, and the numerical approximation.

The problem of approximation of the predictor mapping can be important even in linear systems. For example, in large scale systems, it is difficult to compute the matrix exponential as well as the convolution integrals that involve the matrix exponential. Moreover, when measurements are not available continuously then non-standard control approaches must be used. Indeed, we show again that the problem is solved by means of a hybrid feedback law and the main result for linear systems is given next.

**Theorem 1.2:** *Consider the linear system:*

$$\dot{x}(t) = Ax(t) + Bu(t-\tau)$$
$$x(t) \in \Re^n, u(t) \in \Re^m \quad (1.7)$$

*where $\tau > 0$, $A \in \Re^{n \times n}$, $B \in \Re^{n \times m}$. Let $k \in \Re^{m \times n}$ be such that $(A+Bk)$ is Hurwitz. For every $r > 0$ there exists an integer $N^*(r) \geq 1$ such that for every $N \geq N^*(r)$ there exist constants $Q, \sigma > 0$ with the following property: for every partition $\{T_i\}_{i=0}^{\infty}$ of $\Re_+$ with $\sup_{i \geq 0}(T_{i+1} - T_i) \leq r$, the solution of the closed-loop system (1.7) with*

$$\dot{z}(t) = (A+Bk)z(t), \; z(t) \in \Re^n$$
$$u(t) = kz(t) \quad \text{, for } t \in [T_i, T_{i+1}) \quad (1.8)$$

$$z(T_i) = \left(I + \frac{\tau}{N}A\right)^N x(T_i) + \sum_{k=0}^{N-1}\left(I + \frac{\tau}{N}A\right)^{N-1-k} B \int_{k\frac{\tau}{N}}^{(k+1)\frac{\tau}{N}} u(T_i - \tau + s)ds \quad (1.9)$$

*and arbitrary initial condition $x(0) = x_0 \in \Re^n$, $u(s) = u_0(s)$ for $s \in [-\tau, 0)$ satisfies the estimate:*

$$|x(t)| + \sup_{t-\tau \leq s < t}|u(s)| \leq Q\exp(-\sigma t)\left(|x_0| + \sup_{-\tau \leq s < 0}|u_0(s)|\right), \; \forall t \geq 0 \quad (1.10)$$

Again, notice that (1.9) is the output $z_N$ of the repeated application of the explicit Euler numerical scheme to control system (1.7), i.e., the application of the algorithm

$$z_{j+1} = (I + hA)z_j + B \int_{jh}^{(j+1)h} u(T_i - \tau + s)ds, \text{ for } j = 0, \ldots, N-1$$

with $z_0 = x(T_i)$ and step size $h = \frac{\tau}{N}$. Theorem 1.2 is proved by means of a small-gain methodology, which is different from the proof of Theorem 1.1, and its proof is constructive. It should be emphasized that Theorem 1.2 cannot be obtained by the specialization of Theorem 1.1 to the linear case. Moreover, in Section 4 the control practitioner can actually find explicit formulae for the computation of $N^*(r) \geq 1$.

The structure of the paper is as follows: Section 2 provides some results for the numerical explicit Euler scheme for control systems, which are necessary for the proofs of the main results. The results in Section 2 are not available in numerical analysis textbooks but their proofs are made in the same way with the corresponding results for systems without inputs. Section 3 is devoted to the proof of Theorem 1.1. Section 4 is devoted to the proof of Theorem 1.2 and the presentation of a simple example. Section 5 provides the concluding remarks of the present work. The Appendix contains the proofs of certain auxiliary results.

*Notation.* Throughout the paper we adopt the following notation:

* For a vector $x \in \Re^n$ we denote by $|x|$ its usual Euclidean norm, by $x'$ its transpose. The norm $|A|$ of a matrix $A \in \Re^{m \times n}$ is defined by $|A| = \max\{|Ax| : x \in \Re^n, |x| = 1\}$. $I \in \Re^{n \times n}$ denotes the unit matrix. For $x \in \Re^n$ and $\varepsilon > 0$ we denote by $B_\varepsilon(x)$ the closed ball or radius $\varepsilon > 0$ centered at $x \in \Re^n$, i.e., $B_\varepsilon(x) := \{y \in \Re^n : |y - x| \leq \varepsilon\}$.

* $\Re_+$ denotes the set of non-negative real numbers. $Z_+$ denotes the set of non-negative integers. For every $t \geq 0$, $[t]$ denotes the integer part of $t \geq 0$, i.e., the largest integer being less or equal to $t \geq 0$. A partition $\pi = \{T_i\}_{i=0}^\infty$ of $\Re^+$ is an increasing sequence of times with $T_0 = 0$ and $T_i \to +\infty$.

* We say that an increasing continuous function $\gamma : \Re^+ \to \Re^+$ is of class $K$ if $\gamma(0) = 0$. We say that an increasing continuous function $\gamma : \Re^+ \to \Re^+$ is of class $K_\infty$ if $\gamma(0) = 0$ and $\lim_{s \to +\infty} \gamma(s) = +\infty$. By $KL$ we denote the set of all continuous functions $\sigma = \sigma(s,t) : \Re^+ \times \Re^+ \to \Re^+$ with the properties: (i) for each $t \geq 0$ the mapping $\sigma(\cdot, t)$ is of class $K$; (ii) for each $s \geq 0$, the mapping $\sigma(s, \cdot)$ is non-increasing with $\lim_{t \to +\infty} \sigma(s,t) = 0$.

* By $C^j(A)$ ($C^j(A;\Omega)$), where $A \subseteq \Re^n$ ($\Omega \subseteq \Re^m$), $j \geq 0$ is a non-negative integer, we denote the class of functions (taking values in $\Omega \subseteq \Re^m$) that have continuous derivatives of order $j$ on $A \subseteq \Re^n$.

* Let $x : [a - r, b) \to \Re^n$ with $b > a \geq 0$ and $r \geq 0$. By $\breve{T}_r(t)x$ we denote the "open history" of $x$ from $t - r$ to $t$, i.e., $(\breve{T}_r(t)x)(\theta) := x(t + \theta); \theta \in [-r, 0)$, for $t \in [a,b)$.

* Let $I \subseteq \Re^+ := [0, +\infty)$ be an interval. By $L^\infty(I;U)$ we denote the space of measurable and bounded functions $u(\cdot)$ defined on $I$ and taking values in $U \subseteq \Re^m$. Notice that we do not identify functions in $L^\infty(I;U)$ which differ on a measure zero set. For $x \in L^\infty([-r,0];\Re^n)$ we define $\|x\|_r := \sup_{\theta \in [-r,0)} |x(\theta)|$. Notice that $\sup_{\theta \in [-r,0]} |x(\theta)|$ is not the essential supremum but the actual supremum and that is why the quantities $\sup_{\theta \in [-r,0]} |x(\theta)|$ and $\sup_{\theta \in [-r,0)} |x(\theta)|$ do not coincide in general.

* The shift operator $\delta_\tau u$ maps each function $u : [-\tau, 0) \to U$ to the function $\delta_\tau u : [0, \tau) \to U$ with $(\delta_\tau u)(s) = u(-\tau + s)$ for all $s \in [-\tau, 0)$.

* A function $f : A \to \Re$, where $0 \in A \subseteq \Re^n$ is positive definite if $f(0) = 0$ and $f(x) > 0$ for all $x \neq 0$. A function $f : \Re^n \to \Re$ is radially unbounded if the set $\{x \in \Re^n : f(x) \leq M\}$ is bounded or empty for every $M > 0$.

## 2. Numerical Approximation of the Solutions of Forward Complete Systems

We consider system (1.1) under the following assumptions:

**(H1)** $f : \Re^n \times \Re^m \to \Re^n$ *is a locally Lipschitz vector field with* $f(0,0) = 0$ *that satisfies:*

$$|f(x,u) - f(y,u)| \leq L(|x| + |y| + |u|)|x - y|, \text{ for all } x, y \in \Re^n, u \in \Re^m \quad (2.1)$$

$$|f(x,u)| \leq (|x| + |u|)L(|x| + |u|), \text{ for all } x \in \Re^n, u \in \Re^m \quad (2.2)$$

*where* $L : \Re_+ \to [1, +\infty)$ *is a continuous, non-decreasing function.*

**(H2)** *System (1.1) is forward complete.*

Assumptions (H1) and (H2) have important consequences for system (1.1). Next we point out two consequences which will be used in this section:

**(C1)** *There exist a $C^2$ function $W: \Re^n \to [1,+\infty)$ which is radially unbounded, a constant $c > 0$ and a function $p \in K_\infty$ such that*

$$\nabla W(x) f(x,u) \leq cW(x) + p(|u|), \text{ for all } x \in \Re^n, u \in \Re^m \quad (2.3)$$

**(C2)** *For every $\tau > 0$ there exists a function $a_\tau \in K_\infty$ such that the solution $x(t)$ of (1.1) with arbitrary initial condition $x(0) = x_0$ corresponding to arbitrary measurable and essentially bounded input $u:[0,\tau) \to \Re^m$ satisfies*

$$|x(t)| \leq a_\tau(|x_0| + \|u\|), \text{ for all } t \in [0,\tau] \quad (2.4)$$

*where*

$$\|u\| := \operatorname{ess\,sup}_{t \in [0,\tau)} |u(t)|$$

*Moreover, for every $\tau > 0$, there exists a constant $M_\tau > 0$ such that $a_\tau(s) = M_\tau s$ for all $s \in [0,1]$.*

The existence of a $C^2$ function $W: \Re^n \to [1,+\infty)$ and the existence of a function $a_\tau \in K_\infty$ satisfying the requirements of assumption (C2) are direct consequences of Theorem 1, Corollary 2.3 in [1] and assumption (H1).

Let $P: \Re_+ \to \Re_+$ be a non-decreasing continuous function that satisfies:

$$P(s) \geq s^2 L^2(s) \max\left\{ |\nabla^2 W(\xi)| : |\xi| \leq s(1 + \tau L(s)) \right\}, \text{ for all } s \geq 0 \quad (2.5)$$

Let $Q: \Re_+ \to \Re_+$ be a non-decreasing continuous function that satisfies:

$$Q(s) \geq 1 + \max\left\{ |x| : W(x) \leq \exp(2c\tau) \max_{|y| \leq s}(W(y)) + \frac{\exp(2c\tau)-1}{2c} p(s) \right\}, \text{ for all } s \geq 0 \quad (2.6)$$

Define for all $s \geq 0$:

$$A(s) := L(Q(s) + a_\tau(s) + s) \quad (2.7)$$

$$B(s) := L(Q(s) + a_\tau(s) + s)(a_\tau(s) + s)L(a_\tau(s) + s) \quad (2.8)$$

Consider the following numerical scheme, which is an extension of the explicit Euler method to systems with inputs: we select a positive integer $N$ and define

$$x_{i+1} = x_i + \int_{ih}^{(i+1)h} f(x_i, u(s))ds, \text{ for } i = 0,\ldots,N-1 \quad (2.9)$$

for $h = \dfrac{\tau}{N}$.

**Theorem 2.1:** *Consider system (1.1) under assumptions (H1), (H2). Let $\tau > 0$ be a positive constant and let a $C^2$ function $W : \Re^n \to [1,+\infty)$ which is radially unbounded, a constant $c > 0$, functions $p \in K_\infty$, $a_\tau \in K_\infty$ be such that assumptions (C1) and (C2) hold. Let $P : \Re_+ \to \Re_+$, $Q : \Re_+ \to \Re_+$, $A : \Re_+ \to \Re_+$, $B : \Re_+ \to \Re_+$ be continuous functions that satisfy (2.5), (2.6), (2.7), (2.8). Let arbitrary $x_0 \in \Re^n$ and arbitrary measurable and essentially bounded input $u : [0,\tau) \to \Re^m$. If*

$$N \geq \tau \frac{P(Q(|x_0| + \|u\|) + \|u\|)}{2c}$$

*then the following inequalities hold:*

$$|x(\tau) - x_N| \leq \frac{\tau B(|x_0| + \|u\|)}{2NA(|x_0| + \|u\|)} \left(\exp(\tau A(|x_0| + \|u\|)) - 1\right) \quad (2.10)$$

$$|x_i| \leq Q(|x_0| + \|u\|), \text{ for all } i = 0,1,\ldots,N \quad (2.11)$$

*where $x(\tau)$ is the solution of (1.1) with initial condition $x(0) = x_0$ corresponding to input $u : [0,\tau) \to \Re^m$ at time $t = \tau$.*

**Remark 2.2:** Inequality (2.10) shows that if we know the initial condition $x(0) = x_0$ and the applied input $u : [0,\tau) \to \Re^m$ then we can estimate all quantities involved in (2.10). Moreover, if we want the approximation error to be less than $\varepsilon > 0$ it suffices to select the positive integer $N$ so that:

$$N \geq \tau \max\left(\frac{B(|x_0| + \|u\|)}{2\varepsilon A(|x_0| + \|u\|)} \left(\exp(\tau A(|x_0| + \|u\|)) - 1\right), \frac{P(Q(|x_0| + \|u\|) + \|u\|)}{2c}\right)$$

Notice that the right hand-side of the above inequality can be evaluated before we start applying the scheme (2.9). The restriction is imposed in order to obtain the uniform bound provided by (2.11) and it is necessary for the control of the increase of the function $W$ (exactly in the same spirit as step size control was applied in [7] for the control of the decrease of the Lyapunov function). The bound provided by (2.11) is useful for the proof of Theorem 1.1.

The proof of Theorem 2.1 depends on three technical lemmas which are stated below and are proved at the Appendix.

**Lemma 2.3:** *Consider system (1.1) under the assumptions of Theorem 2.1. If $|x_i| + \|u\| > 0$ and $h \leq \frac{2cW(x_i)}{P(|x_i| + \|u\|)}$, where $P : \Re_+ \to \Re_+$ is the function involved in (2.5), then*

$$W(x_{i+1}) \leq \exp(2ch)W(x_i) + \int_{ih}^{(i+1)h} \exp(2c(ih + h - s))p(|u(s)|)ds \quad (2.12)$$

**Lemma 2.4:** *Consider system (1.1) under the assumptions of Theorem 2.1. If $h \leq \frac{2c}{P(Q(|x_0| + \|u\|) + \|u\|)}$ then*

$$W(x_i) \leq \exp(2cih)W(x_0) + \int_0^{ih} \exp(2c(ih - s))p(|u(s)|)ds \text{ for all } i = 0,\ldots,N \quad (2.13)$$

*where $Q : \Re_+ \to \Re_+$ is the function involved in (2.6).*

**Lemma 2.5:** *Consider system (1.1) under the assumptions of Theorem 2.1. Define $e_i := x_i - x(ih)$, $i \in \{0,...,N\}$, where $x(t)$ is the solution of (1.1) with initial condition $x(0) = x_0$ corresponding to input $u:[0,\tau) \to \Re^m$ and suppose that $h \leq \dfrac{2c}{P(Q(|x_0|+\|u\|)+\|u\|)}$. Then*

$$|e_i| \leq \frac{h^2}{2} B(|x_0|+\|u\|) \frac{\exp(ihA(|x_0|+\|u\|))-1}{\exp(hA(|x_0|+\|u\|))-1}, \text{ for all } i \in \{1,...,N\} \quad (2.14)$$

where the functions $A, B : \Re_+ \to \Re_+$ are defined by (2.7), (2.8).

We are now ready to provide the proof of Theorem 2.1.

**Proof of Theorem 2.1:** All assumptions of Lemma 2.4 and Lemma 2.5 hold. Consequently, inequalities (2.13), (2.14) hold. Inequality (2.10) follows from using the fact $\exp(hA(|x_0|+\|u\|))-1 \geq hA(|x_0|+\|u\|)$ and definition $h = \dfrac{\tau}{N}$ in conjunction with (2.14) for $i = N$. Moreover, inequality (2.13) implies $W(x_i) \leq \exp(2c\tau)W(x_0) + \dfrac{\exp(2c\tau)-1}{2c}p(\|u\|)$. The previous inequality in conjunction with (2.6) implies (2.11). The proof is complete. ◁

Theorem 2.1 allows us to construct mappings which approximate the solution of (1.1) $\tau$ time units ahead with guaranteed accuracy level. Indeed, let $R \in C^0(\Re_+;\Re_+)$ be a positive definite function with $\liminf\limits_{s \to 0^+} \dfrac{R(s)}{s} > 0$. Define the mapping $\Phi : \Re^n \times L^\infty([0,\tau);\Re^m) \to \Re^n$ by means of the equation:

$$\Phi(x_0, u) := x_N \quad (2.15)$$

where $x_i$, $i = 1,...,N$ are defined by the numerical scheme (2.9) with $h = \dfrac{\tau}{N}$ and

$$N(x_0, \|u\|) := \left\lfloor \tau \max\left( \frac{a_\tau(|x_0|+\|u\|)+\|u\|}{2R(|x_0|+\|u\|)} L(a_\tau(|x_0|+\|u\|)+\|u\|)(\exp(\tau A(|x_0|+\|u\|))-1), \frac{P(Q(|x_0|+\|u\|)+\|u\|)}{2c} \right) \right\rfloor + 1 \quad (2.16)$$

for $|x_0|+\|u\| > 0$ and

$$N(0,0) := 1 \quad (2.17)$$

for $|x_0| = \|u\| = 0$.

By virtue of (2.10) the mapping $\Phi : \Re^n \times L^\infty([0,\tau);\Re^m) \to \Re^n$ satisfies

$$|\Phi(x_0, u) - x(\tau)| \leq R(|x_0|+\|u\|) \quad (2.18)$$

Inequalities (2.10), (2.11) in conjunction with (2.18) and (2.4) implies the following inequality:

$$|\Phi(x_0, u)| \leq \min\left( R(|x_0|+\|u\|) + a_\tau(|x_0|+\|u\|), Q(|x_0|+\|u\|) \right) \quad (2.19)$$

Notice that the mapping $N(x_0, \|u\|)$ defined by (2.16) and (2.17) is locally bounded. Indeed, there exists a constant $M_\tau > 0$ such that $a_\tau(|x_0| + \|u\|) = M_\tau |x_0| + M_\tau \|u\|$ for all $x_0 \in \Re^n$ and for every measurable and essentially bounded input $u : [0, \tau) \to \Re^m$ with $|x_0| + \|u\|$ sufficiently small. Therefore, continuity of all functions involved in (2.16) in conjunction with the fact that $\liminf_{s \to 0^+} \frac{R(s)}{s} > 0$ implies that

$$\sup_{|x_0| + \|u\| \leq s} N(x_0, \|u\|) < +\infty, \text{ for all } s \geq 0 \tag{2.20}$$

Therefore, we conclude:

**Corollary 2.6:** *Consider system (1.1) under the assumptions of Theorem 2.1. For every positive definite function $R \in C^0(\Re_+; \Re_+)$ with $\liminf_{s \to 0^+} \frac{R(s)}{s} > 0$ and for every $\tau > 0$, consider the mapping $\Phi : \Re^n \times L^\infty([0, \tau); \Re^m) \to \Re^n$ defined by (2.15) for all $(x_0, u) \in \Re^n \times L^\infty([0, \tau); \Re^m)$, where $x_i$, $i = 1, \ldots, N$ are defined by the numerical scheme (2.9) with $h = \frac{\tau}{N}$ and $N := N(x_0, \|u\|)$ is defined by (2.16), (2.17). Then inequalities (2.18), (2.19) hold for all $(x_0, u) \in \Re^n \times L^\infty([0, \tau); \Re^m)$, where $x(t)$ denotes the solution of (1.1) with initial condition $x(0) = x_0$ corresponding to input $u : [0, \tau) \to \Re^m$ and $\|u\| := \operatorname*{ess\,sup}_{t \in [0, \tau)} |u(t)|$. Moreover, inequality (2.20) holds for all $s \geq 0$.*

## 3. Proof of Theorem 1.1

This section is devoted to the proof of Theorem 1.1. The proof of Theorem 1.1 is constructive and formulae will be given next for the locally bounded mapping $N : \Re^n \times \Re_+ \to \{1,2,3,\ldots\}$ involved in the hybrid dynamic feedback law defined by (1.3), (1.4) and (1.5). In order to simplify the procedure of the proof we break the proof up into two steps.

First Step: Construction of feedback

Second Step: Rest of proof

The control practitioner, who is not interested in reading the details of the proof, may read only the first step of the proof.

First Step: Construction of feedback

The feedback law is entirely given by (1.3)-(1.5), except for the function $N : \Re^n \times \Re_+ \to \{1,2,3,\ldots\}$, whose construction we give here. We assume the knowledge of a function $L : \Re_+ \to [1, +\infty)$, a $C^2$ function $W : \Re^n \to [1, +\infty)$ and a function $a_\tau \in K_\infty$ satisfying the requirements of assumptions (C1), (C2) of Section 2. As remarked in the previous section, the existence of a function $L : \Re_+ \to [1, +\infty)$, a $C^2$ function $W : \Re^n \to [1, +\infty)$ and a function $a_\tau \in K_\infty$ satisfying the requirements of assumptions (H1), (C1), (C2) are direct consequences of Theorem 1, Corollary 2.3 in [1] and the fact that $f : \Re^n \times \Re^m \to \Re^n$ is a continuously differentiable mapping with $f(0,0) = 0$.

Moreover, we need to assume the knowledge of a Lyapunov function for the closed-loop system (1.1) with $u(t) = k(x(t))$. More specifically, we assume the existence of a positive definite, radially unbounded function $V \in C^1(\Re^n; \Re^+)$ for (2.1), constants $\varepsilon, K, \mu > 0$ and a function $\rho \in K_\infty$ such that the following hold:

$$\nabla V(x) f(x, k(x)) \leq -\rho(V(x)), \quad \forall x \in \Re^n \tag{3.1}$$

$$|x|^2 \leq V(x) \leq K|x|^2, \quad \forall x \in B_\varepsilon(0) \tag{3.2}$$

$$|\nabla V(x)| \leq 2K|x|, \quad \forall x \in B_\varepsilon(0) \tag{3.3}$$

$$\nabla V(x) f(x, k(x)) \leq -\mu|x|^2, \quad \forall x \in B_\varepsilon(0) \tag{3.4}$$

The existence of a Lyapunov function for the closed-loop system (1.1) with $u(t) = k(x(t))$ satisfying (3.1), (3.2), (3.3), (3.4) is a direct consequence of Proposition 4.4 in [7].

Based on the knowledge of all the functions and constants described above, we next proceed to the construction of new functions. The first functions to define are the continuous, non-decreasing functions $P: \Re_+ \to \Re_+$, $Q: \Re_+ \to \Re_+$, $A: \Re_+ \to \Re_+$, $B: \Re_+ \to \Re_+$ that satisfy (2.5), (2.6), (2.7), (2.8). Next, we define:

- functions $a_i \in K_\infty$ ($i = 1,...,4$) and constants $k_1, k_2, k_3, k_4 > 0$ that satisfy:

$$a_1(|x|) \leq V(x) \leq a_2(|x|), \text{ for all } x \in \Re^n \tag{3.5}$$

$$|\nabla V(x)| \leq a_3(|x|) \text{ and } |k(x)| \leq a_4(|x|), \text{ for all } x \in \Re^n \tag{3.6}$$

$$a_1(s) := k_1 s^2, \; a_2(s) := k_2 s^2, \; a_3(s) := k_3 s, \; a_4(s) := k_4 s, \text{ for all } s \in [0, \varepsilon] \tag{3.7}$$

- a continuous, non-decreasing function $M: \Re_+ \to [1, +\infty)$ that satisfies:

$$|f(x, k(z)) - f(x, k(x))| \leq M(|x| + |z|)|z - x|, \text{ for all } z, x \in \Re^n \tag{3.8}$$

The reader should notice that the existence of functions $M: \Re_+ \to [1, +\infty)$, $a_i \in K_\infty$ ($i = 1,...,4$) and constants $k_1, k_2, k_3, k_4 > 0$ satisfying (3.5), (3.6), (3.7) and (3.8) is a direct consequence of (a) the fact that $V \in C^1(\Re^n; \Re^+)$ is positive definite and radially unbounded (see Lemma 3.5 in [10]), (b) of Lemma 2.4 in [6], (c) of inequalities (3.2), (3.3) and (d) of the fact that $f: \Re^n \times \Re^m \to \Re^n$ and $k: \Re^n \to \Re^m$ are continuously differentiable mappings with $k(0) = 0$.

Moreover, define for all $s \geq 0$:

$$D_r(s) := a_3(a_r(s) + s) M(a_r(s) + s) \exp(rL(a_r(s) + s)), \; q(s) := a_4(a_1^{-1}(a_2(s))) + a_1^{-1}(a_2(s)) \tag{3.9}$$

where $a_r \in K_\infty$ is the function involved in (2.4) with $\tau$ replaced by $r > 0$ and $L: \Re_+ \to [1, +\infty)$ is the function involved in assumption (H1).

Next select a constant $\delta > 0$, such that:

$$a_1^{-1}(a_2(2a_1^{-1}(\delta))) \leq \varepsilon \tag{3.10}$$

Having selected $\delta > 0$, we are in a position to select a constant $\gamma > 0$, so that:

$$\gamma \leq \min\left(a_1^{-1}(\delta), \frac{1}{2}\rho\left(\frac{\delta}{2}\right)\right) \tag{3.11}$$

Define:

$$\phi := k_3 M\left(a_1^{-1}\left(a_2\left(2a_1^{-1}(\delta)\right)\right) + a_1^{-1}(\delta)\right)\exp(r\tilde{L}) \tag{3.12}$$

$$\tilde{L} := L\left((1+k_4)a_1^{-1}\left(a_2\left(2a_1^{-1}(\delta)\right)\right) + a_1^{-1}(\delta)\right) \tag{3.13}$$

and moreover, select a constant $\tilde{R} > 0$, so that:

$$k_4\sqrt{\frac{k_2}{k_1}}\tilde{R} < 1, \quad \frac{\sqrt{2}k_2\phi\left(\sqrt{k_1} + k_4\sqrt{k_2}\right) + \mu k_1 k_4\sqrt{k_2}}{\mu k_1\sqrt{k_1}}\tilde{R} < 1 \tag{3.14}$$

Finally, define:

$$R(s) := \min\left(\frac{\gamma}{\max(1, D_r(a_\tau(s) + q(Q(s))))}, \tilde{R}s, \frac{1}{2}a_2^{-1}\left(a_1\left(a_4^{-1}\left(\frac{s}{2}\right)\right)\right)\right) \tag{3.15}$$

Notice that by virtue of (3.7) and the fact that $Q(s) \geq 1$ for all $s \geq 0$, it follows from definition (3.15) that $\liminf_{s \to 0^+} \frac{R(s)}{s} = \min\left(\tilde{R}, \frac{1}{4k_4}\sqrt{\frac{k_1}{k_2}}\right) > 0$. Therefore, Corollary 2.6 guarantees that the mapping $N : \Re^n \times \Re_+ \to \{1,2,3,...\}$ defined by (2.16), (2.17) is locally bounded and the mapping $\Phi : \Re^n \times L^\infty([0,\tau); \Re^m) \to \Re^n$ defined by (2.15) satisfies inequalities (2.18), (2.19) for all $(x_0, u) \in \Re^n \times L^\infty([0,\tau); \Re^m)$, where $x(t)$ denotes the solution of (1.1) with initial condition $x(0) = x_0$ corresponding to input $u : [0,\tau) \to \Re^m$ and $\|u\| := \underset{t \in [0,\tau)}{\operatorname{ess\,sup}}|u(t)|$.

Second Step: Rest of proof

Having completed the design of the feedback law by constructing the function $N : \Re^n \times \Re_+ \to \{1,2,3,...\}$ in (2.16), we are now ready to prove some basic results concerning the closed-loop system (1.2) with (1.3), (1.4), (1.5) and (2.16).

The following claim shows that practical stabilization is achieved. Its proof is provided in the Appendix.

**Claim 1:** *There exists $\sigma \in KL$ such that for every partition $\{T_i\}_{i=0}^\infty$ of $\Re_+$ with $\sup_{i \geq 0}(T_{i+1} - T_i) \leq r$, for every $x_0 \in \Re^n$ and $u_0 \in L^\infty([-\tau, 0); \Re^m)$, the solution of (1.2), (1.3), (1.4) and (1.5) with initial condition $x(0) = x_0$, $\tilde{T}_\tau(0)u = u_0$ satisfies the following inequality for all $t \geq 0$:*

$$V(x(t)) \leq \max\left(\sigma\left(|x_0| + \|u_0\|_\tau, t\right), \rho^{-1}(2\gamma)\right) \tag{3.16}$$

*where $\rho \in K_\infty$ is the function involved in (3.1) and $\gamma > 0$ is the constant involved in (3.11) and (3.15).*

The following claim shows that local exponential stabilization is achieved. Its proof is provided in the Appendix.

**Claim 2:** *There exist constants $Q_1, Q_2, \omega > 0$ such that for each partition $\{T_i\}_{i=0}^{\infty}$ of $\Re_+$ with $\sup_{i \geq 0}(T_{i+1} - T_i) \leq r$, for each $x_0 \in \Re^n$ and $u_0 \in L^{\infty}([-\tau,0);\Re^m)$, the solution of (1.2), (1.3), (1.4) and (1.5) with initial condition $x(0) = x_0$, $\breve{T}_\tau(0)u = u_0$ satisfies the following inequalities:*

$$|u(t)|\exp(\omega(t - T_j)) \leq Q_1\left(\sup_{T_j \leq w \leq T_j + \tau}(|x(w)|) + \|\breve{T}_\tau(T_j)u\|_\tau\right), \text{ for all } t \geq T_j \qquad (3.17)$$

$$|x(t)|\exp(\omega(t - T_j - \tau)) \leq Q_2\left(\sup_{T_j \leq w \leq T_j + \tau}(|x(w)|) + \|\breve{T}_\tau(T_j)u\|_\tau\right), \text{ for all } t \geq T_j + \tau \qquad (3.18)$$

*where $T_j$ is the smallest sampling time for which it holds $V(x(T_j + \tau)) \leq \delta$, where $\delta > 0$ is the constant involved in (3.10) and (3.11).*

The following claim guarantees that $u$ is bounded. Its proof is provided in the Appendix.

**Claim 3:** *There exists a non-decreasing function $G: \Re_+ \to \Re_+$ such that for each partition $\{T_i\}_{i=0}^{\infty}$ of $\Re_+$ with $\sup_{i \geq 0}(T_{i+1} - T_i) \leq r$, for each $x_0 \in \Re^n$ and $u_0 \in L^{\infty}([-\tau,0);\Re^m)$, the solution of (1.2), (1.3), (1.4) and (1.5) with initial condition $x(0) = x_0$, $\breve{T}_\tau(0)u = u_0$ satisfies the following inequality for all $t \geq 0$:*

$$|x(t)| + \|\breve{T}_\tau(t)u\|_\tau \leq G(|x_0| + \|u_0\|_\tau) \qquad (3.19)$$

We are now ready to prove Theorem 1.1. Let arbitrary partition $\{T_i\}_{i=0}^{\infty}$ of $\Re_+$ with $\sup_{i \geq 0}(T_{i+1} - T_i) \leq r$, $x_0 \in \Re^n$, $u_0 \in L^{\infty}([-\tau,0);\Re^m)$ and consider the solution of (1.2), (1.3), (1.4) and (1.5) with (arbitrary) initial condition $x(0) = x_0$, $\breve{T}_\tau(0)u = u_0$.

Inequalities (3.5) and (2.4) imply that the smallest sampling time $T_j$ for which $V(x(T_j + \tau)) \leq \delta$ holds is $T_0 = 0$ for the case $a_2(a_\tau(|x_0| + \|u_0\|_\tau)) \leq \delta$. Moreover, the fact that there exists a constant $M_\tau > 0$ such that $a_\tau(s) = M_\tau s$ for all $s \in [0,1]$, in conjunction with inequalities (3.17), (3.18), allow us to conclude that that there exists a constant $\Omega > 0$

$$|x(t)| + \|\breve{T}_\tau(t)u\|_\tau \leq \Omega\exp(-\omega t)(|x_0| + \|u_0\|_\tau), \text{ for all } t \geq 0 \qquad (3.20)$$

provided that $|x_0| + \|u_0\|_\tau \leq \min\left(1, \frac{1}{M_\tau}a_2^{-1}(\delta)\right)$.

Proposition 7 in [25] in conjunction with (3.16), (3.11) and the fact that $\sup_{i \geq 0}(T_{i+1} - T_i) \leq r$, allow us to guarantee the existence of a non-decreasing function $\widetilde{T}: \Re_+ \to \Re_+$ such that the smallest sampling time $T_j$ for which $V(x(T_j + \tau)) \leq \delta$ holds satisfies $T_j \leq \widetilde{T}(|x_0| + \|u_0\|_\tau)$ for all $(x_0, u) \in \mathfrak{R}^n \times L^\infty([0,\tau);\mathfrak{R}^m)$. Combining (3.17), (3.18), (3.19) with the previous inequality, allows us to conclude the existence of a non-decreasing function $\widetilde{G}: \mathfrak{R}_+ \to \mathfrak{R}_+$ such that the following inequality holds for all $t \geq 0$:

$$|x(t)| + \|\breve{T}_\tau(t)u\|_\tau \leq \exp(-\omega t)\widetilde{G}(|x_0| + \|u_0\|_\tau) \qquad (3.21)$$

Consequently, using (3.20) and (3.21) we conclude that (1.6) holds with $C(s) := \frac{1}{s}\int_s^{2s} \widetilde{C}(w)dw$ for all $s > 0$ and $C(0) := 0$, where

$$\widetilde{C}(s) := \max\left(1, \frac{\widetilde{G}(l)}{\Omega l}\right)\Omega s, \text{ for all } s \in [0,l]$$

$$\widetilde{C}(s) := \max\left(\Omega s, \widetilde{G}(s)\right), \text{ for all } s > l$$

$$l := \min\left(1, \frac{1}{M_\tau}a_2^{-1}(\delta)\right)$$

The proof of Theorem 1.1 is complete. ◁

## 4. Linear Systems

We first start this section by providing the formulae and some explanations for Theorem 1.2. In this section we prove that the integer $N^* = N^*(r)$ in the statement of Theorem 1.2 can be selected as the smallest integer that satisfies the inequality

$$\tau|k|\exp(|A|r)(|A|\gamma\exp(|A|\tau) + |B|(|A|\tau\exp(|A|\tau) + 1)(1 + \gamma|k|))(\exp(|A|\tau) - 1) < 2N^* \qquad (4.1)$$

where $\gamma > 0$ is a constant for which the estimate

$$|x(t)| \leq M\exp(-\lambda t)|x_0| + \gamma \sup_{0 \leq s \leq t}\left(\exp(-\lambda(t-s))|u(s-\tau) - kx(s)|\right) \qquad (4.2)$$

holds for the solution $x(t)$ of $\dot{x}(t) = (A + Bk)x(t)$ with arbitrary initial condition $x(0) = x_0 \in \mathfrak{R}^n$ corresponding to arbitrary measurable and locally essentially bounded input $u:[-\tau,+\infty) \to \mathfrak{R}^m$, for certain appropriate constants $M, \lambda > 0$. The constant $\gamma > 0$ can be selected to be any constant with $\gamma > \frac{\sqrt{|B'PB|}}{\mu}$, where $P \in \mathfrak{R}^{n \times n}$ is a symmetric positive definite matrix and $\mu > 0$ is a constant that satisfies $P(A+Bk) + (A+Bk)'P + 2\mu P \leq 0$ and $P \geq I$.

It is clear that the term $\left(I+\frac{\tau}{N}A\right)^N x(T_i)$ in (1.9) is an approximation of $\exp(A\tau)x(T_i)$. Moreover, the term $\sum_{k=0}^{N-1}\left(I+\frac{\tau}{N}A\right)^{N-1-k} B \int_{k\frac{\tau}{N}}^{(k+1)\frac{\tau}{N}} u(T_i-\tau+s)ds$ in (1.9) is an approximation of

$$\int_0^\tau \exp(A(\tau-s))Bu(T_i-\tau+s)dw = \sum_{k=0}^{N-1} \int_{k\frac{\tau}{N}}^{(k+1)\frac{\tau}{N}} \exp(A(\tau-s))Bu(T_i-\tau+s)dw.$$

Therefore $z(T_i)$ is an approximation of $x(T_i+\tau) = \exp(A\tau)x(T_i) + \int_0^\tau \exp(A(\tau-s))Bu(T_i-\tau+s)dw$. The formula (1.9) is similar to the formula (11.7) on page 177 of the book [29]. However, there is a difference between formula (1.9) and formula (11.7) on page 177 of [29]: in formula (1.9) the matrix exponentials $\exp(A\tau)$ and $\exp\left(A\left(\tau-(k+1)\frac{\tau}{N}\right)\right)$ are approximated by $\left(I+\frac{\tau}{N}A\right)^N$ and $\left(I+\frac{\tau}{N}A\right)^{N-1-k}$, respectively. In other words, formula (1.9) does not require the computation of the matrix exponentials. This feature can be important for large-scale systems.

Another advantage of our analysis is the computation of explicit bounds for the integer $N^* = N^*(r)$. In [29], it is shown that the error converges to zero when $N \to +\infty$ (Theorem 11.6 on page 187) but no explicit bound is provided.

Finally, it must be noted that the controller (1.8), (1.9) requires measurements at discrete time instants (the sampling times) and does not require continuous measurement of the state vector. Moreover, it guarantees robustness to perturbations of the sampling schedule: estimate (1.10) holds for all partitions $\{T_i\}_{i=0}^\infty$ of $\Re_+$ with $\sup_{i\geq 0}(T_{i+1}-T_i) \leq r$.

From this point on, we start proving Theorem 1.2. We first prove the following technical result.

**Lemma 4.1:** *For every integer $N \geq 1$, $x(0) \in \Re^n$ and $u \in L^\infty([0,\tau);\Re^m)$, the solution $x(t)$ of $\dot{x}(t) = Ax(t) + Bu(t)$ satisfies*

$$|x_N - x(\tau)| \leq \frac{h}{2}a(\exp(a\tau)-1)\exp(a\tau)|x(0)| + \frac{h}{2}b(a\tau\exp(a\tau)+1)(\exp(a\tau)-1)\|u\| \qquad (4.3)$$

*where* $x_{j+1} = (I+hA)x_j + B\int_{jh}^{(j+1)h} u(s)ds$ $(j=0,...,N-1)$, $x_0 = x(0) \in \Re^n$, $\|u\| = \sup_{0\leq s<\tau}|u(s)|$, $h = \frac{\tau}{N}$, $a = |A|$ *and* $b = |B|$.

**Proof:** If we define $e_j = x_j - x(jh)$ for $j=0,...,N$ then it follows that

$e_{j+1} = (I+hA)e_j + \int_{jh}^{(j+1)h} A(x(jh)-x(s))ds$ for $j=0,...,N-1$. Therefore, we get

$$|e_{j+1}| \leq |I+hA||e_j| + |A|\int_{jh}^{(j+1)h}|x(jh)-x(s)|ds, \text{ for } j=0,...,N-1 \qquad (4.4)$$

For all $s \in [jh, \tau]$ it holds that $|x(s) - x(jh)| \leq |A| \int_{jh}^{s} |x(w) - x(jh)| dw + (s - jh)|A||x(jh)| + (s - jh)|B|\|u\|$. The previous inequality in conjunction with the Gronwall-Bellman Lemma implies that $|x(s) - x(jh)| \leq (s - jh)(|A||x(jh)| + |B|\|u\|) \exp(|A|(s - jh))$ for all $s \in [jh, \tau]$. The previous inequality in conjunction with (4.4) and the fact that $|I + hA| \leq 1 + h|A|$ gives:

$$|e_{j+1}| \leq (1 + h|A|)|e_j| + \frac{h^2}{2}|A|(|A||x(jh)| + |B|\|u\|)\exp(h|A|), \text{ for } j = 0, \ldots, N-1 \qquad (4.5)$$

Using (4.5) and the facts that $e_0 = 0$ and $\max_{j=0,\ldots,N-1} |x(jh)| \leq \max_{0 \leq s \leq \tau - h} |x(s)|$, we obtain:

$$|e_N| \leq \frac{h}{2}\left(|A| \max_{0 \leq s \leq \tau - h} |x(s)| + |B|\|u\|\right)\exp(h|A|)\left((1 + h|A|)^N - 1\right) \qquad (4.6)$$

For all $s \in [0, \tau]$ it holds that $|x(s)| \leq |x(0)| + |A| \int_0^s |x(w)| dw + s|B|\|u\|$. The previous inequality in conjunction with the Gronwall-Bellman Lemma implies that $|x(s)| \leq (|x(0)| + s|B|\|u\|)\exp(|A|s)$ for all $s \in [0, \tau]$. Consequently, the previous estimate in conjunction with (4.6) gives:

$$|e_N| \leq \frac{h}{2}\left(a|x(0)|\exp(a\tau) + (\tau - h)ab\|u\|\exp(a\tau) + b\|u\|\exp(ah)\right)\left((1 + ha)^N - 1\right) \qquad (4.7)$$

where $a = |A|$ and $b = |B|$. Inequality (4.7) in conjunction with the facts that $1 + ha \leq \exp(ha)$ and $h = \frac{\tau}{N}$ gives us the desired inequality (4.3). ◁

We are now ready to give the proof of Theorem 1.2.

**Proof of Theorem 1.2:** The proof is established by means of a small-gain argument. Let $\sigma \in (0, \lambda]$ be sufficiently small so that

$$|k|\exp((a + \sigma)r + \sigma\tau)(C_2 + (C_2|k| + C_1)\gamma) < 1 \qquad (4.8)$$

where $\gamma, \lambda > 0$ are the constants involved in (4.2), $C_1 = \frac{h}{2}a(\exp(a\tau) - 1)\exp(a\tau)$, $C_2 = \frac{h}{2}b(a\tau \exp(a\tau) + 1)(\exp(a\tau) - 1)$, $a = |A|$ and $b = |B|$. The existence of sufficiently small $\sigma \in (0, \lambda]$ is guaranteed by (4.1).

Let $\{T_i\}_{i=0}^{\infty}$ be an arbitrary partition of $\Re_+$ with $\sup_{i \geq 0}(T_{i+1} - T_i) \leq r$ and consider the solution of the closed-loop system (1.7) with (1.8), (1.9) and arbitrary initial condition $x(0) = x_0 \in \Re^n$, $\breve{T}_\tau(0)u = u_0 \in L^{\infty}([-\tau, 0); \Re^m)$. Notice that inequality (4.2) implies that for all $t \geq 0$ we have:

$$|x(t + \tau)|\exp(\sigma(t + \tau)) \leq M|x(\tau)|\exp(\sigma\tau) + \gamma\exp(\sigma\tau)\sup_{0 \leq w \leq t}(\exp(\sigma w)|u(w) - kx(w + \tau)|) \qquad (4.9)$$

Since $u(t) = kz(t)$ for all $t \geq 0$ we get from (4.9):

$$|x(t + \tau)|\exp(\sigma(t + \tau)) \leq M|x(\tau)|\exp(\sigma\tau) + \gamma|k|\exp(\sigma\tau)\sup_{0 \leq w \leq t}(\exp(\sigma w)|z(w) - x(w + \tau)|) \qquad (4.10)$$

For all $t \in [T_i, T_{i+1})$ it holds that $|z(t) - x(t+\tau)| \leq \exp(ar)|z(T_i) - x(T_i + \tau)|$. Moreover, inequality (4.3) implies that $|z(T_i) - x(\tau + T_i)| \leq C_1|x(T_i)| + C_2 \sup_{T_i - \tau \leq s < T_i} |u(s)|$, where $C_1 = \frac{h}{2}a(\exp(a\tau) - 1)\exp(a\tau)$ and $C_2 = \frac{h}{2}b(a\tau \exp(a\tau) + 1)(\exp(a\tau) - 1)$. The above inequalities in conjunction with the fact that $T_{i+1} \leq T_i + r$ imply that the following inequality holds for all $t \in [T_i, T_{i+1})$:

$$|z(t) - x(t+\tau)|\exp(\sigma t) \leq C_1 \exp((a+\sigma)r)|x(T_i)|\exp(\sigma T_i)$$
$$+ C_2 \exp(ar + \sigma t) \sup_{T_i - \tau \leq s < T_i} |u(s)| \quad (4.11)$$

Notice that since $u(t) = kz(t)$ for all $t \geq 0$, we get $\sup_{T_i - \tau \leq s < T_i} |u(s)| \leq \|u_0\|_\tau + |k| \sup_{0 \leq s \leq T_i} |z(s)|$ for all $i \in Z^+$ with $T_i < \tau$ and $\sup_{T_i - \tau \leq s < T_i} |u(s)| \leq |k| \sup_{T_i - \tau \leq s \leq T_i} |z(s)|$ for all $i \in Z_+$ with $T_i \geq \tau$. Consequently, since $T_{i+1} \leq T_i + r$, we have

$$\exp(\sigma t) \sup_{T_i - \tau \leq s < T_i} |u(s)| \leq \exp(\sigma(T_i + r)) \sup_{T_i - \tau \leq s < T_i} |u(s)| \leq \exp(\sigma(\tau + r))\|u_0\|_\tau + |k|\exp(\sigma(\tau + r)) \sup_{0 \leq s \leq T_i} (\exp(\sigma s)|z(s)|),$$

for all $i \in Z_+$ with $T_i < \tau$ and $t \in [T_i, T_{i+1})$

$$\exp(\sigma t) \sup_{T_i - \tau \leq s < T_i} |u(s)| \leq |k|\exp(\sigma(T_i + r)) \sup_{T_i - \tau \leq s \leq T_i} |z(s)| \leq |k|\exp(\sigma(\tau + r)) \sup_{T_i - \tau \leq s \leq T_i} (\exp(\sigma s)|z(s)|)$$

for all $i \in Z_+$ with $T_i \geq \tau$ and $t \in [T_i, T_{i+1})$

Combining the two above cases and (4.11) it follows that the following estimate holds for all $t \in [T_i, T_{i+1})$:

$$|z(t) - x(t+\tau)|\exp(\sigma t) \leq C_1 \exp((a+\sigma)r)|x(T_i)|\exp(\sigma T_i)$$
$$+ C_2 \exp(ar + \sigma(r+\tau))\|u_0\|_\tau$$
$$+ C_2|k|\exp(ar + \sigma(r+\tau)) \sup_{0 \leq s \leq t} (\exp(\sigma s)|z(s)|) \quad (4.12)$$

The trivial inequalities $\sup_{0 \leq s \leq t}(\exp(\sigma s)|z(s)|) \leq \sup_{0 \leq s \leq t}(\exp(\sigma s)|z(s) - x(s+\tau)|) + \exp(-\sigma\tau) \sup_{0 \leq s \leq t}(\exp(\sigma(s+\tau))|x(s+\tau)|)$ and $|x(T_i)|\exp(\sigma T_i) \leq \exp(\sigma\tau) \max_{0 \leq s \leq \tau}(|x(s)|) + \sup_{0 \leq s \leq t}(\exp(\sigma(s+\tau))|x(s+\tau)|)$ (which holds for all $t \in [T_i, T_{i+1})$) in conjunction with (4.12) imply that the following estimate holds for all $t \geq 0$:

$$\sup_{0 \leq s \leq t}(\exp(\sigma s)|z(s) - x(s+\tau)|) \leq C_1 \exp((a+\sigma)r)\exp(\sigma\tau) \max_{0 \leq s \leq \tau}|x(s)|$$
$$+ C_2 \exp(ar + \sigma(r+\tau))\|u_0\|_\tau$$
$$+ C_2|k|\exp(ar + \sigma(r+\tau)) \sup_{0 \leq s \leq t}(\exp(\sigma s)|z(s) - x(s+\tau)|) \quad (4.13)$$
$$+ (C_2|k| + C_1)\exp((a+\sigma)r) \sup_{0 \leq s \leq t}(\exp(\sigma(s+\tau))|x(s+\tau)|)$$

For all $s \in [0, \tau]$ it holds that $|x(s)| \leq |x(0)| + |A|\int_0^s |x(w)|dw + s|B|\|u_0\|_\tau$. The previous inequality in conjunction with the Gronwall-Bellman Lemma and definitions $a = |A|$ and $b = |B|$ implies that $|x(s)| \leq (|x(0)| + sb\|u_0\|_\tau)\exp(as)$ for all $s \in [0, \tau]$. Consequently, it follows that $\max_{0 \leq s \leq \tau}|x(s)| \leq (|x(0)| + \tau b\|u_0\|_\tau)\exp(a\tau)$. The previous inequality combined with (4.10) and (4.13) implies that the following estimates hold for all $t \geq 0$:

$$\sup_{0 \leq s \leq t} \left( \exp(\sigma s) |z(s) - x(s+\tau)| \right) \leq C_1 \exp((a+\sigma)r) \exp((a+\sigma)\tau) |x(0)|$$
$$+ \exp((a+\sigma)r + \sigma\tau)(C_1 \tau b \exp(a\tau) + C_2) \|u_0\|_\tau$$
$$+ C_2 |k| \exp(ar + \sigma(r+\tau)) \sup_{0 \leq s \leq t} \left( \exp(\sigma s) |z(s) - x(s+\tau)| \right) \quad (4.14)$$
$$+ (C_2 |k| + C_1) \exp((a+\sigma)r) \sup_{0 \leq s \leq t} \left( \exp(\sigma(s+\tau)) |x(s+\tau)| \right)$$

$$\sup_{0 \leq s \leq t} \left( |x(s+\tau)| \exp(\sigma(s+\tau)) \right) \leq M \tau b \exp((a+\sigma)\tau) \|u_0\|_\tau + M |x(0)| \exp((a+\sigma)\tau)$$
$$+ \gamma |k| \exp(\sigma\tau) \sup_{0 \leq w \leq t} \left( \exp(\sigma w) |z(w) - x(w+\tau)| \right) \quad (4.15)$$

Combining (4.14), (4.15) and (4.8) we obtain for all $t \geq 0$:

$$\sup_{0 \leq s \leq t} \left( \exp(\sigma s) |z(s) - x(s+\tau)| \right) \leq \frac{\exp((a+\sigma)(r+\tau))(C_1 + (C_2|k| + C_1)M)}{1 - |k| \exp((a+\sigma)r + \sigma\tau)(C_2 + (C_2|k| + C_1)\gamma)} |x(0)|$$
$$+ \frac{\exp((a+\sigma)r + \sigma\tau)((C_1 \tau b \exp(a\tau) + C_2) + (C_2|k| + C_1) \exp(a\tau) M \tau b)}{1 - |k| \exp((a+\sigma)r + \sigma\tau)(C_2 + (C_2|k| + C_1)\gamma)} \|u_0\|_\tau \quad (4.16)$$

$$\sup_{0 \leq s \leq t} \left( \exp(\sigma(s+\tau)) |x(s+\tau)| \right) \leq$$
$$\left( M \exp((a+\sigma)\tau) + \gamma |k| \frac{\exp((a+\sigma)(r+\tau) + \sigma\tau)(C_1 + (C_2|k| + C_1)M)}{1 - |k| \exp((a+\sigma)r + \sigma\tau)(C_2 + (C_2|k| + C_1)\gamma)} \right) |x(0)| \quad (4.17)$$
$$\left( M \tau b \exp((a+\sigma)\tau) + \gamma |k| \frac{\exp((a+\sigma)r + 2\sigma\tau)((C_1 \tau b \exp(a\tau) + C_2) + (C_2|k| + C_1) \exp(a\tau) M \tau b)}{1 - |k| \exp((a+\sigma)r + \sigma\tau)(C_2 + (C_2|k| + C_1)\gamma)} \right) \|u_0\|_\tau$$

Since $\max_{0 \leq s \leq \tau} |x(s)| \leq (|x(0)| + \tau b \|u_0\|_\tau) \exp(a\tau)$, it follows from (4.17) that:

$$\sup_{0 \leq s \leq t} \left( \exp(\sigma s) |x(s)| \right) \leq$$
$$\left( M \exp((a+\sigma)\tau) + \gamma |k| \frac{\exp((a+\sigma)(r+\tau) + \sigma\tau)(C_1 + (C_2|k| + C_1)M)}{1 - |k| \exp((a+\sigma)r + \sigma\tau)(C_2 + (C_2|k| + C_1)\gamma)} \right) |x(0)| \quad (4.18)$$
$$\left( M \tau b \exp((a+\sigma)\tau) + \gamma |k| \frac{\exp((a+\sigma)r + 2\sigma\tau)((C_1 \tau b \exp(a\tau) + C_2) + (C_2|k| + C_1) \exp(a\tau) M \tau b)}{1 - |k| \exp((a+\sigma)r + \sigma\tau)(C_2 + (C_2|k| + C_1)\gamma)} \right) \|u_0\|_\tau$$

The existence of a constant $Q > 0$ for which (1.10) holds is a direct consequence of (4.16), (4.17), (4.18) and the following inequalities:

$$\exp(\sigma t) \|\tilde{T}_\tau(t) u\|_\tau \leq \exp(\sigma t) \sup_{t-\tau \leq s < t} |u(s)| \leq \exp(\sigma \tau) \|u_0\|_\tau + |k| \exp(\sigma t) \sup_{\max(0, t-\tau) \leq s < t} |z(s)|$$
$$\leq \exp(\sigma \tau) \|u_0\|_\tau + |k| \exp(\sigma t) \sup_{\max(0, t-\tau) \leq s \leq t} |z(s) - x(s+\tau)| + |k| \exp(\sigma t) \sup_{\max(0, t-\tau) \leq s \leq t} |x(s+\tau)|$$
$$\leq \exp(\sigma \tau) \|u_0\|_\tau + |k| \exp(\sigma \tau) \sup_{0 \leq s \leq t} \left( \exp(\sigma s) |z(s) - x(s+\tau)| \right) + |k| \sup_{0 \leq s \leq t} \left( \exp(\sigma(s+\tau)) |x(s+\tau)| \right)$$

which hold for all $t \geq 0$. The proof is complete. ◁

Next we present a simple example, which shows the usefulness of the formulae that were provided in this section.

**Example 4.2:** Consider the scalar system

$$\dot{x}(t) = x(t) + u(t-1), \text{ with } x(t) \in \Re, u(t) \in \Re \tag{4.19}$$

This is the example considered on page 188 of the book [29]. Here $a = 1 = b$ and $k = -p$, where $p > 1$. The constant $\gamma > 0$ can be selected to be any constant with $\gamma > \frac{1}{p-1}$ by using the matrix $P = [1]$. Therefore, inequality (4.1) shows that the integer $N^* = N^*(r)$ can be selected so that

$$N^* > f(p) = \frac{p \exp(r)(e + (e+1)(2p-1))(e-1)}{2(p-1)}$$

Figure 1 shows the behaviour of $f(p)$ with respect to $p > 1$ for $r = 1$. The minimum is attained at $p = 1.93$ with $f(p) \approx 64.71$, showing that at least 65 grid points should be used for the Euler approximation of the solution of the open-loop system over a time window of length 1 with an accuracy guarantee in the form of (4.3). It follows that for every $N \geq N^*(r)$ there exist constants $Q, \sigma > 0$ with the following property:

For every partition $\{T_i\}_{i=0}^{\infty}$ of $\Re_+$ with $\sup_{i \geq 0}(T_{i+1} - T_i) \leq r$, the solution of the closed-loop system (4.19) with

$$\begin{matrix} \dot{z}(t) = (p-1)z(t), z(t) \in \Re \\ u(t) = -pz(t) \end{matrix}, \text{ for } t \in [T_i, T_{i+1}) \tag{4.20}$$

$$z(T_i) = \left(1 + \frac{1}{N}\right)^N x(T_i) + \sum_{k=0}^{N-1}\left(1 + \frac{1}{N}\right)^{N-1-k} \int_{\frac{k}{N}}^{\frac{k+1}{N}} u(T_i - 1 + s)ds \tag{4.21}$$

and arbitrary initial condition $x(0) = x_0 \in \Re^n$, $\breve{T}_\tau(0)u = u_0 \in L^\infty([-1,0); \Re^m)$ satisfies estimate (1.10).

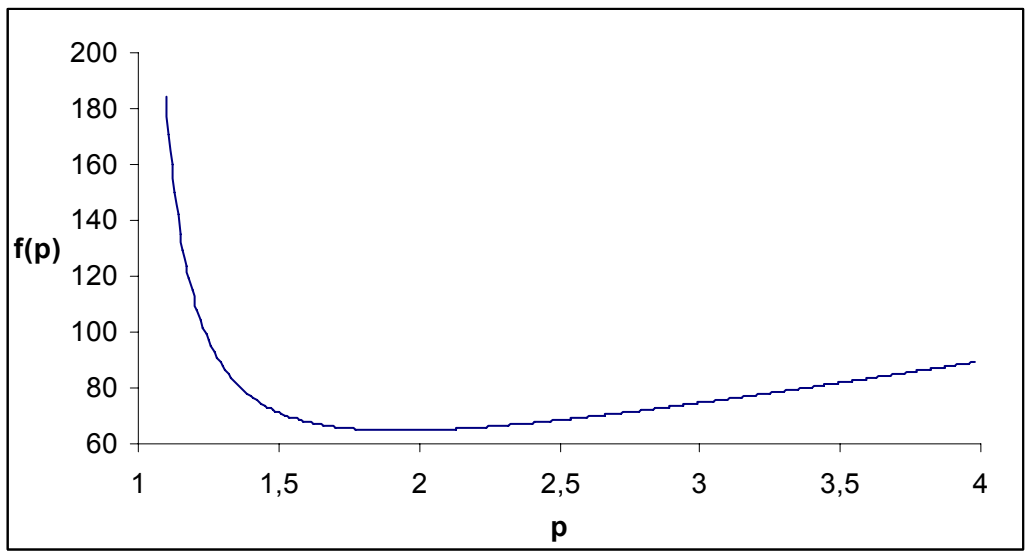

Figure 1: The behaviour of $f(p)$ with respect to $p > 1$ for $r = 1$

# 5. Concluding Remarks

This work has focused on a key aspect of the implementation problem for predictor-based feedback laws: the problem of the approximation of the predictor mapping. It was shown that the numerical approximation of the predictor mapping by means of the explicit Euler numerical scheme in conjunction with a hybrid feedback law that uses sampled measurements, can be used for the global stabilization of all forward complete nonlinear systems that are globally asymptotically stabilizable and locally exponentially stabilizable in the delay-free case.

Special results were provided for the linear time invariant case. Both for the linear and nonlinear case, explicit formulae are provided for the estimates for the design parameters of the resulting hybrid control scheme.

The present paper goes beyond the approximation results in [5] by removing the global Lipschitz restriction.

More remains to be done for the approximations of the integrals involved in the explicit Euler scheme by easily implementable formulae. Results for linear systems are already given in [15,16,17,18,19,20,22,24,27,28,29]. Furthermore, one cannot ignore the possibility of using different numerical schemes (except the explicit Euler scheme; see [2]): the use of implicit numerical schemes may require fewer grid points than the grid points needed for the explicit Euler scheme. Finally, there is the challenging problem of using numerical approximations for cases where the measured output is not necessarily the state vector and there is measurement delay (see [8,9]).

# Appendix

**Proof of Lemma 2.3:** Define the function:

$$g(\lambda) = W(x_i + \lambda(x_{i+1} - x_i)) \tag{A.1}$$

for $\lambda \in [0,1]$. The following equalities hold for all $\lambda \in [0,1]$:

$$\frac{dg}{d\lambda}(\lambda) = \nabla W(x_i + \lambda(x_{i+1} - x_i))(x_{i+1} - x_i)$$

$$\frac{d^2 g}{d\lambda^2}(\lambda) = (x_{i+1} - x_i)' \nabla^2 W(x_i + \lambda(x_{i+1} - x_i))(x_{i+1} - x_i) \tag{A.2}$$

Moreover, notice that by virtue of (2.2) and (2.9), it holds that $|x_{i+1} - x_i| \leq h(|x_i| + \|u\|)L(|x_i| + \|u\|)$. The previous inequality in conjunction with (2.5) and (A.2) gives:

$$\left|\frac{d^2 g}{d\lambda^2}(\lambda)\right| \leq h^2 P(|x_i| + \|u\|) \tag{A.3}$$

where $P: \Re_+ \to \Re_+$ is the function involved in (2.5). Furthermore, inequality (2.3) in conjunction with (2.9) and (A.2) gives:

$$\frac{dg}{d\lambda}(0) = \nabla W(x_i) \int_{ih}^{(i+1)h} f(x_i, u(s))ds \leq chW(x_i) + \int_{ih}^{(i+1)h} p(|u(s)|)ds \tag{A.4}$$

Combining (A.1), (A.3) and (A.4), we get:

$$W(x_{i+1}) = g(1) \leq (1 + ch)W(x_i) + \int_{ih}^{(i+1)h} p(|u(s)|)ds + \frac{h^2}{2} P(|x_i| + \|u\|) \tag{A.5}$$

Inequality (A.5) in conjunction with the following inequality

$$(1 + ch)W(x_i) + \int_{ih}^{(i+1)h} p(|u(s)|)ds + \frac{h^2}{2} P(|x_i| + \|u\|) \leq \exp(2ch)W(x_i) + \int_{ih}^{(i+1)h} \exp(2c(ih + h - s))p(|u(s)|)ds$$

which holds for all $h \leq \frac{2cW(x_i)}{P(|x_i| + \|u\|)}$ imply that (2.12) holds. The proof is complete. ◁

**Proof of Lemma 2.4:** We will first prove that if there exists $j \in \{0, ..., N-1\}$ such that $\|u\| + \min_{i=0,...,j} |x_i| > 0$ and $h \leq \frac{2c}{P(Q(|x_0| + \|u\|) + \|u\|)}$ then (2.13) holds for all $i = 0, ..., j+1$. The proof is by induction.

First notice that (2.13) holds for $i = 0$. Suppose that it holds for some $i \in \{0, ..., j\}$. Clearly, inequality (2.13) implies $W(x_i) \leq \exp(2c\tau)W(x_0) + \frac{\exp(2c\tau) - 1}{2c} p(\|u\|)$. The previous inequality in conjunction with (2.6) implies $|x_i| \leq Q(|x_0| + \|u\|)$. Consequently, the facts that $P: \Re_+ \to \Re_+$ is non-decreasing and $W(x_i) \geq 1$ imply $h \leq \frac{2c}{P(Q(|x_0|+\|u\|)+\|u\|)} \leq \frac{2cW(x_i)}{P(|x_i|+\|u\|)}$. Since $|x_i|+\|u\| > 0$ and $h \leq \frac{2cW(x_i)}{P(|x_i|+\|u\|)}$, Lemma 2.3 shows that:

$$W(x_{i+1}) \leq \exp(2ch)W(x_i) + \int_{ih}^{(i+1)h} \exp(2c(ih+h-s))p(|u(s)|)ds$$

The above inequality in conjunction with (2.13) shows that (2.13) holds for $i$ replaced by $i+1$.

The case that there exists $j \in \{0,...,N-1\}$ with $\|u\| + \min_{i=0,...,j}|x_i| = 0$ can be treated in the following way. Let $j \in \{0,...,N-1\}$ be the smallest integer with $\|u\| + \min_{i=0,...,j}|x_i| = 0$. This implies that $\|u\| = 0$ and $|x_j| = 0$. Since $f(0,0) = 0$, (2.9) implies that $|x_i| = 0$ for all $i = j+1,...,N$. Consequently, (2.13) holds for all $i = j+1,...,N$.

The proof is complete. ◁

**Proof of Lemma 2.5:** Notice that, by virtue of (2.9), the following equation holds for all $i \in \{0,...,N-1\}$:

$$e_{i+1} = e_i + \int_{ih}^{(i+1)h} (f(x_i,u(s)) - f(x(s),u(s)))ds \qquad (A.6)$$

Inequality (2.1) implies the following inequality for all $i \in \{0,...,N-1\}$ and $s \in [ih,(i+1)h]$:

$$|f(x_i,u(s)) - f(x(s),u(s))| \leq L(|x_i|+|x(s)|+\|u\|)|x_i - x(s)| \qquad (A.7)$$

Using the definition $e_i := x_i - x(ih)$ and inequalities (2.2), (2.4) we get for all $i \in \{0,...,N-1\}$ and $s \in [ih,(i+1)h]$:

$$\begin{aligned}|x_i - x(s)| &\leq |e_i| + |x(s) - x(ih)| \\ &\leq |e_i| + (s-ih)(|x(s)|+\|u\|)L(|x(s)|+\|u\|) \\ &\leq |e_i| + (s-ih)(a_\tau(|x_0|+\|u\|)+\|u\|)L(a_\tau(|x_0|+\|u\|)+\|u\|)\end{aligned} \qquad (A.8)$$

Notice that all hypotheses of Lemma 2.4 hold. Therefore inequality (2.13) holds for all $i = 0,...,N$. Clearly, inequality (2.13) implies $W(x_i) \leq \exp(2c\tau)W(x_0) + \frac{\exp(2c\tau)-1}{2c}p(\|u\|)$. The previous inequality in conjunction with (2.6) implies $|x_i| \leq Q(|x_0|+\|u\|)$ for all $i = 0,...,N$. Exploiting the fact that $|x_i| \leq Q(|x_0|+\|u\|)$ for all $i = 0,...,N$ and (A.6), (A.7), (A.8), we obtain for all $i \in \{0,...,N-1\}$:

$$\begin{aligned}|e_{i+1}| &\leq |e_i| + hL(Q(|x_0|+\|u\|)+a_\tau(|x_0|+\|u\|)+\|u\|)|e_i| \\ &+ \frac{h^2}{2}L(Q(|x_0|+\|u\|)+a_\tau(|x_0|+\|u\|)+\|u\|)(a_\tau(|x_0|+\|u\|)+\|u\|)L(a_\tau(|x_0|+\|u\|)+\|u\|)\end{aligned} \qquad (A.9)$$

Definitions (2.7), (2.8) in conjunction with inequality (A.9) shows that the following recursive relation holds for all $i \in \{0,...,N-1\}$

$$|e_{i+1}| \leq \exp(hA(|x_0|+\|u\|))|e_i| + \frac{h^2}{2}B(|x_0|+\|u\|) \tag{A.10}$$

Using the fact $e_0 = 0$, in conjunction with relation (A.10), gives the desired inequality (2.14). The proof is complete. ◁

**Proof of Claim 1:** First we show that for each partition $\{T_i\}_{i=0}^{\infty}$ of $\Re_+$ with $\sup_{i \geq 0}(T_{i+1}-T_i) \leq r$, for each $x_0 \in \Re^n$ and $u_0 \in L^\infty([-\tau,0);\Re^m)$, the solution of (1.2), (1.3), (1.4) and (1.5) with initial condition $x(0)=x_0$, $\breve{T}_\tau(0)u = u_0$ is unique and exists for all $t \geq 0$.

The solution of (1.2), (1.3), (1.4) and (1.5) is determined by the following process:

Initial Step: Given $x(0)=x_0$, $\breve{T}_\tau(0)u = u_0$ we determine the solution $x(t)$ of (1.2) for $t \in [0,\tau]$. Notice that the solution is unique. Inequality (2.4) implies the following estimate:

$$|x(t)| \leq a_\tau(|x_0| + \|u_0\|_\tau), \text{ for all } t \in [0,\tau] \tag{A.11}$$

$i$-th Step: Given $x(t)$ for $t \in [0,T_i+\tau]$ and $u(t)$ for $t \in [-\tau,T_i)$ we determine $x(t)$ for $t \in [0,T_{i+1}+\tau]$ and $u(t)$ for $t \in [-\tau,T_{i+1})$. The solution $z(t)$ of (1.3) for $t \in [T_i,T_{i+1})$ with initial condition $z(T_i)=z_N$ is unique (by virtue of the fact that $f$ and $k$ are locally Lipschitz mappings). Inequality (3.1) implies:

$$V(z(t)) \leq V(z(T_i)), \text{ for all } t \in [T_i,T_{i+1}) \tag{A.12}$$

We determine $u(t)$ for $t \in [T_i,T_{i+1})$ by means of the equation $u(t) = k(z(t))$. Notice that inequalities (3.5), (3.6) in conjunction with (A.12) imply the following inequality for all $t \in [T_i,T_{i+1})$:

$$|u(t)| = |k(z(t))| \leq a_4(a_1^{-1}(a_2(|z(T_i)|))) \tag{A.13}$$

Finally, we determine the solution $x(t)$ of (1.2) for $t \in [T_i+\tau, T_{i+1}+\tau]$. Notice that the solution is unique. The fact that $T_{i+1} - T_i \leq r$ in conjunction with inequality (2.4) with $\tau$ replaced by $r > 0$ and inequality (A.13) implies the estimate:

$$|x(t)| \leq a_r(|x(T_i+\tau)| + a_4(a_1^{-1}(a_2(|z(T_i)|)))), \text{ for all } t \in [T_i+\tau, T_{i+1}+\tau] \tag{A.14}$$

Next we evaluate the difference $z(t) - x(t+\tau)$ for $t \in [T_i, T_{i+1})$. Exploiting (2.1) we get:

$$|z(t) - x(t+\tau)| = \left| z(T_i) - x(T_i+\tau) + \int_{T_i}^t (f(z(s),k(z(s))) - f(x(s+\tau),k(z(s))))ds \right|$$

$$\leq |z(T_i) - x(T_i+\tau)| + \int_{T_i}^t L(|z(s)| + |x(s+\tau)| + |k(z(s))|)|z(s) - x(s+\tau)|ds$$

Using the right inequality (3.5), inequalities (A.12), (A.13), (A.14), in conjunction with the above inequality, we obtain:

$$|z(t) - x(t+\tau)| \leq |z(T_i) - x(T_i + \tau)|$$
$$+ L\left(a_1^{-1}\left(a_2\left(|z(T_i)|\right)\right) + a_4\left(a_1^{-1}\left(a_2\left(|z(T_i)|\right)\right)\right) + a_r\left(|x(T_i+\tau)| + a_4\left(a_1^{-1}\left(a_2\left(|z(T_i)|\right)\right)\right)\right)\right) \int_{T_i}^{t} |z(s) - x(s+\tau)| ds$$

Define $\varphi(s) := a_r(s) + s$. Using the Growall-Bellman lemma, the above inequality and the fact that $T_{i+1} - T_i \leq r$, we get for all $t \in [T_i, T_{i+1})$:

$$|z(t) - x(t+\tau)| \leq |z(T_i) - x(T_i+\tau)| \exp\left(rL\left(\varphi\left(|x(T_i+\tau)| + q\left(|z(T_i)|\right)\right)\right)\right) \quad (A.15)$$

Next we evaluate the quantity $\nabla V(x(t+\tau))f(x(t+\tau), k(z(t)))$ for $t \in [T_i, T_{i+1})$. Using inequality (3.1) we get:

$$\nabla V(x(t+\tau))f(x(t+\tau), k(z(t))) \leq -\rho(V(x(t+\tau))) +$$
$$\nabla V(x(t+\tau))\left(f(x(t+\tau), k(z(t))) - f(x(t+\tau), k(x(t+\tau)))\right)$$

The following estimate follows from (3.6), (3.8) and the above inequality:

$$\nabla V(x(t+\tau))f(x(t+\tau), k(z(t))) \leq -\rho(V(x(t+\tau))) +$$
$$a_3\left(|x(t+\tau)|\right) M\left(|x(t+\tau)| + |z(t)|\right) |x(t+\tau) - z(t)|$$

Using the above inequality in conjunction with inequality (3.5), inequalities (A.12), (A.14) and definitions $q(s) := a_4\left(a_1^{-1}(a_2(s))\right) + a_1^{-1}(a_2(s))$, $\varphi(s) := a_r(s) + s$, we get:

$$\nabla V(x(t+\tau))f(x(t+\tau), k(z(t))) \leq -\rho(V(x(t+\tau))) +$$
$$a_3\left(\varphi\left(|x(T_i+\tau)| + q\left(|z(T_i)|\right)\right)\right) M\left(\varphi\left(|x(T_i+\tau)| + q\left(|z(T_i)|\right)\right)\right) |x(t+\tau) - z(t)| \quad (A.16)$$

Combining inequalities (A.15), (A.16) and definition (3.9) we obtain for all $t \in [T_i, T_{i+1})$:

$$\nabla V(x(t+\tau))f(x(t+\tau), k(z(t))) \leq -\rho(V(x(t+\tau))) + D_r\left(|x(T_i+\tau)| + q\left(|z(T_i)|\right)\right) |x(T_i+\tau) - z(T_i)| \quad (A.17)$$

$$D_r(s) := a_3(a_r(s) + s) M(a_r(s) + s) \exp(rL(a_r(s) + s)) \quad R(s) := \min\left(\frac{\gamma}{D_r(a_\tau(s) + q(Q(s)))}, \tilde{R}s, \frac{1}{2} a_2^{-1}\left(a_1\left(a_4^{-1}\left(\frac{s}{2}\right)\right)\right)\right)$$

Since $z(T_i) = z_N$ (recall (1.4)), it follows from (2.18), (2.19) and (1.2), (1.4) that the following inequalities hold for all $i = 0,1,2,...$:

$$|z(T_i) - x(T_i + \tau)| \leq R\left(|x(T_i)| + \|\breve{T}_\tau(T_i)u\|_\tau\right) \quad (A.18)$$

$$|z(T_i)| \leq Q\left(|x(T_i)| + \|\breve{T}_\tau(T_i)u\|_\tau\right) \quad (A.19)$$

Since $|x(T_i + \tau)| \leq a_\tau\left(|x(T_i)| + \|\breve{T}_\tau(T_i)u\|_\tau\right)$ (recall (2.4)), we obtain from (A.17), (A.18), (A.19) and definition (3.15) for all $t \in [T_i, T_{i+1})$:

$$\frac{d}{dt} V(x(t+\tau)) \leq -\rho(V(x(t+\tau))) + \gamma \quad (A.20)$$

Using (A.20) and Lemma 2.14, page 82 in [6], we obtain for all $t \geq 0$:

$$V(x(t+\tau)) \leq \max\left(\tilde{\sigma}(V(x(\tau)), t), \rho^{-1}(2\gamma)\right) \quad (A.21)$$

for certain function $\tilde{\sigma} \in KL$. Combining (3.5), (A.11) and (A.21) we obtain inequality (3.16) with $\sigma(s,t) := \tilde{\sigma}(a_2(a_\tau(s)), t-\tau)$ for all $t > \tau$ and $\sigma(s,t) := \tilde{\sigma}(a_2(a_\tau(s)), 0)$ for all $t \in [0, \tau]$. The proof is complete. ◁

**Proof of Claim 2:** Let arbitrary partition $\{T_i\}_{i=0}^{\infty}$ of $\mathfrak{R}_+$ with $\sup_{i \geq 0}(T_{i+1} - T_i) \leq r$, $x_0 \in \mathfrak{R}^n$, $u_0 \in L^\infty([-\tau, 0); \mathfrak{R}^m)$ and consider the solution of (1.2), (1.3), (1.4) and (1.5) with (arbitrary) initial condition $x(0) = x_0$, $\tilde{T}_\tau(0)u = u_0$. Inequalities (3.11) and (3.16) guarantee that there exists a unique smallest sampling time $T_j$ such that $V(x(T_j + \tau)) \leq \delta$.

Moreover, inequalities (A.20), (3.11) and (3.5) allow us to conclude that

$$|x(t)| \leq a_1^{-1}(\delta) \text{ and } V(x(t)) \leq \delta, \text{ for all } t \geq T_j + \tau \tag{A.22}$$

Using (A.18), definition (3.15), (3.11) and (A.22) we obtain for all $i \geq j$:

$$|z(T_i)| \leq |z(T_i) - x(T_i + \tau)| + |x(T_i + \tau)| \leq \gamma + a_1^{-1}(\delta) \leq 2a_1^{-1}(\delta) \tag{A.23}$$

Using (A.12), (3.5) and (A.23), we get for all $t \geq T_j$:

$$|z(t)| \leq a_1^{-1}(a_2(2a_1^{-1}(\delta))) \tag{A.24}$$

Next we evaluate the difference $z(t) - x(t + \tau)$ for $t \geq T_j$. Exploiting (2.1) and inequalities (3.6), (3.7), (3.10), (A.22), (A.24) and definition (3.13), we get for all $i \geq j$ and $t \in [T_i, T_{i+1})$:

$$|z(t) - x(t+\tau)| = \left| z(T_i) - x(T_i + \tau) + \int_{T_i}^{t} (f(z(s), k(z(s))) - f(x(s+\tau), k(z(s)))) ds \right|$$

$$\leq |z(T_i) - x(T_i + \tau)| + \tilde{L} \int_{T_i}^{t} |z(s) - x(s+\tau)| ds$$

Using the Growall-Bellman lemma, the above inequality and the fact that $T_{i+1} - T_i \leq r$ imply for all $i \geq j$ and $t \in [T_i, T_{i+1})$:

$$|z(t) - x(t+\tau)| \leq |z(T_i) - x(T_i + \tau)| \exp(r\tilde{L}) \tag{A.25}$$

Next we evaluate the quantity $\nabla V(x(t+\tau)) f(x(t+\tau), k(z(t)))$ for $t \in [T_i, T_{i+1})$. Using inequalities (3.4), (3.5), (3.7), (A.22), (3.10), (3.8), (A.24) and (A.25) and definition (3.12), we get for all $i \geq j$ and $t \in [T_i, T_{i+1})$:

$$\nabla V(x(t+\tau)) f(x(t+\tau), k(z(t))) \leq -\mu k_2^{-1} V(x(t+\tau)) + \phi |x(t+\tau)| |x(T_i + \tau) - z(T_i)| \tag{A.26}$$

Using (3.5), (3.7), (3.10), (A.22) and (A.26) we get for all $i \geq j$ and $t \in [T_i, T_{i+1})$:

$$\dot{V}(t+\tau) \leq -\frac{\mu}{2k_2} V(t+\tau) + \frac{k_2}{2\mu k_1} \phi^2 |x(T_i + \tau) - z(T_i)|^2 \tag{A.27}$$

where $V(t) = V(x(t))$. Using (A.18), (3.15) and (A.27) we get for all $i \geq j$ and $t \in [T_i, T_{i+1})$:

$$\dot{V}(t+\tau) \leq -\frac{\mu}{2k_2}V(t+\tau) + \frac{k_2}{\mu k_1}\phi^2 \widetilde{R}^2 |x(T_i)|^2 + \frac{k_2}{\mu k_1}\phi^2 \widetilde{R}^2 \|\widetilde{T}_\tau(T_i)u\|_\tau^2 \quad \text{(A.28)}$$

Let $\omega < \frac{\mu}{4k_2}$ be a positive constant sufficiently small so that

$$k_4\sqrt{\frac{k_2}{k_1}}\widetilde{R}\exp(\omega(r+\tau)) < 1 \text{ and } \sqrt{2}\frac{k_2}{k_1}\phi\widetilde{R}\frac{\exp(\omega(r+\tau))}{\sqrt{\mu^2-4\omega\mu k_2}}\left(1 + \frac{k_4\sqrt{k_2}\exp(\omega(r+\tau))(\widetilde{R}+\exp(-\omega\tau))}{\sqrt{k_1}-\widetilde{R}k_4\sqrt{k_2}\exp(\omega(r+\tau))}\right) < 1 \quad \text{(A.29)}$$

The existence of $0 < \omega < \frac{\mu}{4k_2}$ satisfying (A.29) is guaranteed by (3.14). Using (A.28) and the fact that $\sup_{i\geq 0}(T_{i+1}-T_i) \leq r$, we obtain for all $i \geq j$ and $t \in [T_i, T_{i+1})$:

$$\dot{V}(t+\tau) \leq -\frac{\mu}{2k_2}V(t+\tau) + \frac{k_2}{\mu k_1}\phi^2\widetilde{R}^2 \exp(-2\omega t)\exp(2\omega r) \sup_{T_i \leq s \leq t}\left(\exp(2\omega s)|x(s)|^2\right)$$
$$+ \frac{k_2}{\mu k_1}\phi^2\widetilde{R}^2 \exp(-2\omega t)\exp(2\omega(r+\tau)) \sup_{T_i-\tau \leq s \leq t}\left(\exp(2\omega s)|u(s)|^2\right) \quad \text{(A.30)}$$

The differential inequality (A.30) allows us to conclude that the following differential inequality holds for $t \geq T_j$ a.e.:

$$\dot{V}(t+\tau) \leq -\frac{\mu}{2k_2}V(t+\tau) + \frac{k_2}{\mu k_1}\phi^2\widetilde{R}^2 \exp(-2\omega t)\exp(2\omega r) \sup_{T_j \leq s \leq t}\left(\exp(2\omega s)|x(s)|^2\right)$$
$$+ \frac{k_2}{\mu k_1}\phi^2\widetilde{R}^2 \exp(-2\omega t)\exp(2\omega(r+\tau)) \sup_{T_j-\tau \leq s \leq t}\left(\exp(2\omega s)|u(s)|^2\right) \quad \text{(A.31)}$$

Integrating (A.31) and since $\omega < \frac{\mu}{4k_2}$, we obtain for all $t \geq T_j$:

$$V(t+\tau) \leq \exp(-2\omega(t-T_j))V(T_j+\tau) + \frac{2k_2^2}{\mu k_1}\phi^2\widetilde{R}^2\frac{\exp(-2\omega t)}{\mu-4\omega k_2}\exp(2\omega r)\sup_{T_j \leq s \leq t}\left(\exp(2\omega s)|x(s)|^2\right)$$
$$+ \frac{2k_2^2}{\mu k_1}\phi^2\widetilde{R}^2\frac{\exp(-2\omega t)}{\mu-4\omega k_2}\exp(2\omega(r+\tau))\sup_{T_j-\tau \leq s \leq t}\left(\exp(2\omega s)|u(s)|^2\right) \quad \text{(A.32)}$$

Using (3.5), (3.7), (A.21) and the fact that $\omega < \frac{\mu}{4k_2}$, we obtain from (A.32) for all $t \geq T_j$:

$$|x(t+\tau)|\exp(\omega(t+\tau)) \leq \exp(\omega(T_j+\tau))\sqrt{\frac{k_2}{k_1}}|x(T_j+\tau)|$$
$$+ \sqrt{2}\frac{k_2}{k_1}\phi\widetilde{R}\frac{\exp(\omega(r+\tau))}{\sqrt{\mu^2-4\omega\mu k_2}}\sup_{T_j \leq s \leq t}\left(\exp(\omega s)|x(s)|\right) \quad \text{(A.33)}$$
$$+ \sqrt{2}\frac{k_2}{k_1}\phi\widetilde{R}\frac{\exp(\omega(r+2\tau))}{\sqrt{\mu^2-4\omega\mu k_2}}\sup_{T_j-\tau \leq s \leq t}\left(\exp(\omega s)|u(s)|\right)$$

Using (3.5), (3.6), (3.7), (3.10), (3.15), (A.12), (A.18) and (A.24) we obtain for all $i \geq j$ and $t \in [T_i, T_{i+1})$:

$$|u(t)| = |k(z(t))| \leq k_4 |z(t)| \leq k_4 \sqrt{\frac{k_2}{k_1}} |z(T_i)|$$

$$\leq k_4 \sqrt{\frac{k_2}{k_1}} |z(T_i) - x(T_i + \tau)| + k_4 \sqrt{\frac{k_2}{k_1}} |x(T_i + \tau)| \quad (A.34)$$

$$\leq \widetilde{R} k_4 \sqrt{\frac{k_2}{k_1}} |x(T_i)| + \widetilde{R} k_4 \sqrt{\frac{k_2}{k_1}} \|\widetilde{T}_\tau(T_i) u\|_\tau + k_4 \sqrt{\frac{k_2}{k_1}} |x(T_i + \tau)|$$

Inequality (A.34) in conjunction with the fact that $\sup_{i \geq 0}(T_{i+1} - T_i) \leq r$ implies:

$$|u(t)| \exp(\omega t) \leq \widetilde{R} k_4 \sqrt{\frac{k_2}{k_1}} \exp(\omega r) |x(T_i)| \exp(\omega T_i)$$

$$+ \widetilde{R} k_4 \sqrt{\frac{k_2}{k_1}} \exp(\omega(r+\tau)) \sup_{T_i - \tau \leq s < T_i} (\exp(\omega s) |u(s)|) + k_4 \sqrt{\frac{k_2}{k_1}} \exp(\omega(r-\tau)) |x(T_i + \tau)| \exp(\omega(T_i + \tau))$$

Therefore, we get from the above inequality for all $t \geq T_j$:

$$|u(t)| \exp(\omega t) \leq k_4 \exp(\omega r) \sqrt{\frac{k_2}{k_1}} (\widetilde{R} + \exp(-\omega \tau)) \sup_{T_j - \tau \leq s \leq t} (\exp(\omega(s+\tau)) |x(s+\tau)|)$$

$$+ \widetilde{R} k_4 \sqrt{\frac{k_2}{k_1}} \exp(\omega(r+\tau)) \sup_{T_j - \tau \leq s \leq t} (\exp(\omega s) |u(s)|) \quad (A.35)$$

Distinguishing the cases $\sup_{T_j - \tau \leq s \leq t}(\exp(\sigma s)|u(s)|) = \sup_{T_j \leq s \leq t}(\exp(\sigma s)|u(s)|)$ and $\sup_{T_j - \tau \leq s \leq t}(\exp(\sigma s)|u(s)|) = \sup_{T_j - \tau \leq s < T_j}(\exp(\sigma s)|u(s)|)$ we obtain from (A.35) for all $t \geq T_j$:

$$|u(t)| \exp(\omega t) \leq \frac{k_4 \exp(\omega r) \sqrt{k_2} (\widetilde{R} + \exp(-\omega \tau))}{\sqrt{k_1} - \widetilde{R} k_4 \sqrt{k_2} \exp(\omega(r+\tau))} \sup_{T_j - \tau \leq s \leq t} (\exp(\omega(s+\tau)) |x(s+\tau)|)$$

$$+ \widetilde{R} k_4 \sqrt{\frac{k_2}{k_1}} \exp(\omega(r+\tau)) \sup_{T_j - \tau \leq s < T_j} (\exp(\omega s) |u(s)|) \quad (A.36)$$

Combining (A.33) and (A.36) we get for all $t \geq T_j$:

$$|x(t+\tau)| \exp(\omega(t+\tau)) \leq \exp(\omega(T_j + \tau)) \sqrt{\frac{k_2}{k_1}} |x(T_j + \tau)|$$

$$+ \sqrt{2} \frac{k_2}{k_1} \phi \widetilde{R} \frac{\exp(\omega(r+\tau))}{\sqrt{\mu^2 - 4\omega\mu k_2}} \left(1 + \frac{k_4 \sqrt{k_2} \exp(\omega(r+\tau))(\widetilde{R} + \exp(-\omega \tau))}{\sqrt{k_1} - \widetilde{R} k_4 \sqrt{k_2} \exp(\omega(r+\tau))}\right) \sup_{T_j - \tau \leq s \leq t} (\exp(\omega(s+\tau)) |x(s+\tau)|)$$

$$+ \sqrt{2} \frac{k_2}{k_1} \phi \widetilde{R} \frac{\exp(\omega(r+2\tau))}{\sqrt{\mu^2 - 4\omega\mu k_2}} \sup_{T_j - \tau \leq s < T_j} (\exp(\omega s) |u(s)|)$$

Distinguishing the cases $\sup_{T_j-\tau \leq s \leq t}(\exp(\sigma(s+\tau))|x(s+\tau)|) = \sup_{T_j-\tau \leq s \leq T_j}(\exp(\sigma(s+\tau))|x(s+\tau)|)$, $\sup_{T_j-\tau \leq s \leq t}(\exp(\sigma(s+\tau))|x(s+\tau)|) = \sup_{T_j \leq s \leq t}(\exp(\sigma(s+\tau))|x(s+\tau)|)$ and using the above inequality, we obtain for all $t \geq T_j$:

$$|x(t+\tau)|\exp(\omega(t+\tau)) \leq \frac{\exp(\omega(T_j+\tau))}{1-A}\sqrt{\frac{k_2}{k_1}}|x(T_j+\tau)| \\ + A\sup_{T_j-\tau \leq s \leq T_j}(\exp(\omega(s+\tau))|x(s+\tau)|) + \sqrt{2}\frac{k_2}{k_1(1-A)}\phi\tilde{R}\frac{\exp(\omega(r+2\tau))}{\sqrt{\mu^2-4\omega\mu k_2}}\sup_{T_j-\tau \leq s < T_j}(\exp(\omega s)|u(s)|) \quad (A.37)$$

where $A = \sqrt{2}\frac{k_2}{k_1}\phi\tilde{R}\frac{\exp(\omega(r+\tau))}{\sqrt{\mu^2-4\omega\mu k_2}}\left(1+\frac{k_4\sqrt{k_2}\exp(\omega(r+\tau))(\tilde{R}+\exp(-\omega\tau))}{\sqrt{k_1}-\tilde{R}k_4\sqrt{k_2}\exp(\omega(r+\tau))}\right)$. Inequalities (A.36), (A.37) imply that there exist constants $Q_1, Q_2 > 0$ such that (3.17), (3.18) hold.

The proof is complete. ◁

**Proof of Claim 3:** Let arbitrary partition $\{T_i\}_{i=0}^{\infty}$ of $\Re_+$ with $\sup_{i \geq 0}(T_{i+1}-T_i) \leq r$, $x_0 \in \Re^n$, $u_0 \in L^{\infty}([-\tau,0);\Re^m)$ and consider the solution of (1.2), (1.3), (1.4) and (1.5) with (arbitrary) initial condition $x(0) = x_0$, $\breve{T}_\tau(0)u = u_0$.

Define:
$$b(s) := a_4(a_1^{-1}(a_2(s))), \text{ for all } s \geq 0 \quad (A.38)$$

and notice that $b \in K_\infty$. Moreover, notice that definitions (A.38) and (3.15) imply that

$$R(s) \leq \frac{1}{2}b^{-1}\left(\frac{s}{2}\right), \text{ for all } s \geq 0 \quad (A.39)$$

Furthermore, definition (A.38) and inequality (A.13) imply the following inequality for all $i \in Z_+$ and $t \in [T_i, T_{i+1})$:

$$|u(t)| \leq b(|z(T_i)|) \quad (A.40)$$

Inequalities (3.5), (3.16) imply the existence of a non-decreasing function $g: \Re_+ \to \Re_+$ such that:

$$|x(t)| \leq g(|x_0|+\|u_0\|_\tau), \text{ for all } t \geq 0 \quad (A.41)$$

By virtue of (A.18), (A.39) and (A.41) we get for all $i \in Z_+$:

$$|z(T_i) - x(T_i+\tau)| \leq R\left(|x(T_i)| + \sup_{T_i-\tau \leq s < T_i}(|u(s)|)\right) \\ \leq \frac{1}{2}b^{-1}\left(\frac{1}{2}|x(T_i)| + \frac{1}{2}\sup_{T_i-\tau \leq s < T_i}(|u(s)|)\right) \\ \leq \max\left(\frac{1}{2}b^{-1}(|x(T_i)|), \frac{1}{2}b^{-1}\left(\sup_{T_i-\tau \leq s < T_i}|u(s)|\right)\right) \\ \leq \max\left(\frac{1}{2}b^{-1}(g(|x_0|+\|u_0\|_\tau)), \frac{1}{2}b^{-1}\left(\sup_{T_i-\tau \leq s < T_i}|u(s)|\right)\right)$$

The above inequality in conjunction with (A.41) gives for all $i \in Z_+$:

$$|z(T_i)| \leq |x(T_i + \tau)| + \max\left(\frac{1}{2}b^{-1}\left(g(|x_0| + \|u_0\|_\tau)\right), \frac{1}{2}b^{-1}\left(\sup_{T_i - \tau \leq s < T_i} |u(s)|\right)\right)$$

$$\leq g(|x_0| + \|u_0\|_\tau) + \max\left(\frac{1}{2}b^{-1}\left(g(|x_0| + \|u_0\|_\tau)\right), \frac{1}{2}b^{-1}\left(\sup_{T_i - \tau \leq s < T_i} |u(s)|\right)\right)$$

$$\leq \max\left(2g(|x_0| + \|u_0\|_\tau), b^{-1}\left(g(|x_0| + \|u_0\|_\tau)\right), b^{-1}\left(\sup_{T_i - \tau \leq s < T_i} |u(s)|\right)\right)$$

$$\leq \max\left(2g(|x_0| + \|u_0\|_\tau), b^{-1}\left(g(|x_0| + \|u_0\|_\tau)\right), b^{-1}\left(\sup_{-\tau \leq s < T_i} |u(s)|\right)\right)$$

Furthermore, using (A.40) and the above inequality, we obtain for all $i \in Z_+$:

$$\sup_{T_i \leq s < T_{i+1}} |u(s)| \leq \max\left(\tilde{g}(|x_0| + \|u_0\|_\tau), \sup_{-\tau \leq s < T_i} |u(s)|\right) \quad \text{(A.42)}$$

where $\tilde{g}(s) := \max(2g(s), b^{-1}(g(s)))$ for all $s \geq 0$, is a non-decreasing function. Define the sequence:

$$F_i := \sup_{-\tau \leq s < T_i} |u(s)| \quad \text{(A.43)}$$

Notice that definition (A.43) and the fact that $\sup_{-\tau \leq s < T_{i+1}} |u(s)| = \max\left(\sup_{T_i \leq s < T_{i+1}} |u(s)|, \sup_{-\tau \leq s < T_i} |u(s)|\right)$ in conjunction with (A.42) imply the following inequality for all $i \in Z_+$:

$$F_{i+1} \leq \max\left(\tilde{g}(|x_0| + \|u_0\|_\tau), F_i\right) \quad \text{(A.44)}$$

Inequality (A.44) in conjunction with the fact that $F_0 := \|u_0\|_\tau$ allow us to prove by induction that the following inequality holds for all $i \in Z_+$:

$$F_i \leq \max\left(\tilde{g}(|x_0| + \|u_0\|_\tau), \|u_0\|_\tau\right) \quad \text{(A.45)}$$

Inequality (A.41) in conjunction with inequality (A.45) and definition (A.43) imply that estimate (3.19) holds with $G(s) := g(s) + \max(\tilde{g}(s), s)$ for all $s \geq 0$.
The proof is complete. ◁